\documentclass[12pt, twoside]{article}
\usepackage{amssymb, amsmath, amsthm, verbatim}
\usepackage[cp1251]{inputenc}
\usepackage[english]{babel}
\newcounter{parag}
\newcommand{\sect}[1]
{\refstepcounter{parag}
\begin{center} { \bf\S\,\theparag. #1} \end{center}}

\newtheorem*{theorem*}{Theorem}
\newtheorem*{conjecture*}{Conjecture}
\newtheorem{lemma}{Lemma}[parag]
\newtheorem{prop}[lemma]{Proposition}

\theoremstyle{definition}
\newtheorem*{prf}{Proof}

\begin{document}

\begin{center} {\large\bf Locally finite groups with bounded centralizer chains}
\end{center}

\begin{center} Alexandr Buturlakin and Andrey V. Vasil$'$ev\footnote{The work was partially supported by RFBR Grant 13-01-00505}
\end{center}

\textbf{Abstract.} The $c$-dimension of a group $G$ is the maximal length of a chain of nested centralizers in $G$. We prove that a locally finite group of finite $c$-dimension $k$ has less than $5k$ nonabelian composition factors.


\begin{center}\textbf{Introduction}\end{center}

Let $G$ be a group and $C_G(X)$ be the centralizer of a subset $X$ of~$G$. Since $C_G(X)<C_G(Y)$ if and only if $C_G(C_G(X))>C_G(C_G(Y))$, it follows that the minimal and the maximal conditions for centralizers are equivalent. Thus the length of every chain of nested centralizers in a group with the minimal condition for centralizers is finite. If a uniform bound for the lengths of chains of centralizers of a group $G$ exists, then we refer to maximal such length as $c$-dimension of $G$ following \cite{MS}. The same notion is also known as the height of the lattice of centralizers. It is worth to observe that the class of groups of finite $c$-dimension includes abelian groups, torsion-free hyperbolic groups, linear groups over fields and so on. In addition, it is closed under taking subgroups and finite direct products, but the $c$-dimension of a homomorphic image of a group from this class is not necessary finite.

In 1979 R.\,Bryant and B.\,Hartley \cite{BH} proved that a periodic locally soluble group with the minimal condition for centralizers is soluble. In 2009 E.\,I.\,Khukhro published the paper \cite{Khukhro}, where, in particular, he proved that a periodic locally soluble group of finite $c$-dimension $k$ has the derived length bounded in terms of~$k$. The same paper contains the conjecture attributed to A.\,V.\,Borovik, which asserts that the number of nonabelian composition factors of a locally finite group of finite $c$-dimension~$k$ is bounded in terms of~$k$. The purpose of our work is to prove this conjecture.

\begin{theorem*} Let $G$ be a locally finite group of $c$-dimension $k$. Then the number of nonabelian composition factors of $G$ is less than $5k$.
\end{theorem*}

\sect{Preliminaries}

Given a locally finite group $G$, denote by $\eta(G)$ the number of nonabelian composition factors of $G$.

The following well-known fact (see, for example, \cite[Corollary 3.5]{Meier}) helps us to derive the theorem from the corresponding statement for finite groups.

\begin{lemma}\label{lfsg} If $G$ is a locally finite locally soluble simple group, then $G$ is cyclic.
\end{lemma}

Recall that the factor group of a finite group $G$ by its soluble radical $R$ is an automorphism group of a direct product of nonabelian simple groups. Thus, if the socle $Soc(G/R)$ is a direct product of nonabelian simple groups $S_1$, $S_2$, $\dots$, $S_n$, then $G/R$ is a subgroup of the semidirect product $\left(\operatorname{\operatorname{Aut}}(S_1)\times \operatorname{Aut}(S_2)\times\dots\times \operatorname{Aut}(S_n)\right)\leftthreetimes Sym_n$, where $Sym_n$ permutes $S_1$, $S_2$, $\dots$, $S_n$. By the classification of the finite simple groups, the group of outer automorphisms of a finite simple group is soluble. Therefore, every nonabelian composition factor of $G$ is either a composition factor of $Soc(G/R)$, or a composition factor of the corresponding subgroup of $Sym_n$.

Next three lemmas give an upper bound for the number of nonabelian composition factors of a subgroup of $Sym_n$. We denote by $\mu(G)$ the degree of the minimal faithful permutation representation of a finite group~$G$.

\begin{lemma}[\cite{Holt}, Theorem 2]\label{DegBySolRad} Let $G$ be a finite group. Let $\mathfrak{L}$ be a class of finite groups closed under taking subgroups, homomorphic images and extensions. If $N$ is the maximal normal $\mathfrak{L}$-subgroup of $G$, then $\mu(G)\geqslant\mu(G/N)$.
\end{lemma}

\begin{lemma}[\cite{Praeger}, Theorem 3.1]\label{DirProdSimple} Let $S_1$, $S_2$, $\dots$, $S_r$ be simple groups. Then $\mu(S_1\times S_2\times\dots\times S_r)=\mu(S_1)+\mu(S_2)+\dots+\mu(S_r)$.
\end{lemma}

\begin{lemma}\label{CompFactorSn} If $G$ is a subgroup of a symmetric group $Sym_n$, then $\eta(G)\leqslant(n-1)/4$.
\end{lemma}

\begin{prf} We proceed by induction on $n$. If $R$ is the soluble radical of $G$, then Lemma~\ref{DegBySolRad} implies that $\mu(G/R)$ does not exceed $\mu(G)$. Hence, we may assume that the soluble radical of $G$ is trivial. Let the socle $Soc(G)$ of $G$ be the direct product of nonabelian simple groups $S_1$, $S_2$, $\dots$, $S_l$. It follows from Lemma~\ref{DirProdSimple} that $l\leqslant n/5$. Again $G$ is a subgroup of the semidirect product $\left(\operatorname{Aut}(S_1)\times \operatorname{Aut}(S_2)\times\dots\times \operatorname{Aut}(S_l)\right)\leftthreetimes Sym_l$. By inductive hypothesis, $\eta(G)\leqslant n/5+(n/5-1)/4=(n-1)/4$.
\end{prf}

\noindent\textsc{Remark.} The group $Sym_{n}$, where $n=5^k$ with $k\geqslant 1$, contains a subgroup $G$ isomorphic to the permutation wreath product $(\dots((Alt_5\wr Alt_5)\wr Alt_5)\dots)$, where the wreath product is applied $k-1$ times. We have $\eta(G)=\frac{5^k-1}{5-1}=\frac{n-1}{4}$.

\medskip

The following lemma is a key for bounding the number of composition factors of $Soc(G/R)$ for a finite group $G$.

\begin{lemma}[\cite{Khukhro}, Lemma 3]\label{LemmaKhukhro} If an elementary abelian $p$-group $E$ of order $p^n$ acts faithfully on a finite nilpotent $p'$-group $Q$, then there exists a series of subgroups $E=E_0>E_1>E_2>\dots>E_n=1$ such that all inclusions  $C_Q(E_0)<C_Q(E_1)<\dots<C_Q(E_n)$ are strict.
\end{lemma}

As usual, $O_p(G)$ stands for the largest normal $p$-subgroup of a finite group $G$, while
$O_{p'}(G)$ denotes the largest normal $p'$-subgroup of $G$. If a series of commutator subgroups of
a group $G$ stabilizes, then we denote by $G^{(\infty)}$ the last subgroup of this series. A
quasisimple group is a perfect central extension of a nonabelian simple group. The layer $E(G)$ is
the subgroup of $G$ generated by all subnormal quasisimple subgroups of $G$, the latter are called
components of $G$. Recall that the layer is a central product of components of $G$.

\newpage

\sect{Proofs}

\begin{prop}\label{FiniteCase} Let $G$ be a finite group of $c$-dimension $k$. Then $\eta(G)<5k$.
\end{prop}

\begin{prf} Let $R$ be the soluble radical of $G$. If $P$ is a Sylow subgroup of $R$, then $G/R \simeq N_G(P)/(R\cap N_G(P))$, so nonabelian composition factors of $N_G(P)$ and $G$ coincide. On the other hand, $c$-dimension of $N_G(P)$ as a subgroup of $G$ is at most $k$. Therefore, we may assume that $N_G(P)=G$ for every Sylow subgroup $P$ of $R$, i.\,e. that $R$ is nilpotent.

Obviously, we suppose that $R\neq G$. Put $\overline{G}=G/R$. The socle $\overline{L}$ of $\overline{G}$ is the direct product of nonabelian simple groups $S_1$, $S_2,\dots$, $S_n$. As observed in preliminaries, the group $\overline{G}/\overline{L}$ is an extension of a normal soluble subgroup by a subgroup of the symmetric group $Sym_n$. By Lemma \ref{CompFactorSn}, an arbitrary subgroup of $Sym_n$ has less than $n/4$ nonabelian composition factors. Thus, it is sufficient to show that $\eta(\overline{L})=n\leqslant 4k$. In particular, we may assume that $G$ coincides with $L$, the preimage of $\overline{L}$ in $G$, and nonabelian composition factor of $G$ are the groups $S_1$, $S_2,\dots$, $S_n$.

Let $K=C_G(R)$. The normal subgroup $\overline{K}=KR/R$ of $\overline{G}$ is a direct product of
nonabelian simple group. Without loss of generality, we may suppose that $\overline{K}=S_1\times
S_2\times\ldots\times S_l$ for some $1\leq l\leq n$. For $i=1,\ldots,l$ denote by $K_i$ the
preimage of $S_i$ in $K$. Then subgroup $H_i=K_i^{(\infty)}$ is normal in $K$ and is a perfect
central extension of $S_i$, so it is a component of $K$. Therefore, if $E(K)$ is the layer of $K$,
then $KR=E(K)R$ and $E(K)$ is a central product of $H_1$, $H_2,\ldots$, $H_l$. Hence
$\eta(K)=\eta(E(K))=l$. Since $[H_i,H_j]=1$ for $i \neq j$, all inclusions
$C_{E(K)}(H_1)<C_{E(K)}(H_1H_2)<\dots<C_{E(K)}(H_1H_2\dots H_l)$ are strict. Thus, $l\leqslant k$.

Let $P$ be a Sylow $p$-subgroup of $G$ and $\overline{P}$ be the image of $P$ in $\overline{G}$.
Since $O_p(R)\leq C_G(O_{p'}(R))$, the action of $P$ on $O_{p'}(R)$ by conjugation induces the
action of $\overline{P}$ on~$O_{p'}(R)$. Given a prime $p$, define the set $\mathcal{F}_p$ as
follows: a subgroup $S_i$ of $\overline{G}$ lies in $\mathcal{F}_p$ whenever there is an element
$g$ of order $p$ in $S_i$ acting faithfully on~$O_{p'}(R)$. Lemma~\ref{LemmaKhukhro} yields that
$|\mathcal{F}_p|\leqslant k$ for every prime $p$. On the other hand, if $S_i$ does not lie in
$\mathcal{F}_p$, then $S_i$ is a subgroup of $C_G(O_{p'}(R))R/R$. It follows from the
classification of finite simple groups that the order of every nonabelian finite simple group is an
even number which is a multiple of $3$ or $5$. Since $R=O_2(R)\times O_{2'}(R)$, every $S_i$ either
belongs to $\mathcal{F}_2\cup\mathcal{F}_3\cup\mathcal{F}_5$, or is a subgroup of
$\overline{K}=C_G(R)R/R$. Thus,
$\eta(G)\leq|\mathcal{F}_2|+|\mathcal{F}_3|+|\mathcal{F}_5|+\eta(K)\leq4k$, as required.
\end{prf}

\textbf{Proof of the theorem.} Now $G$ is locally finite group. Assume $\eta(G)\geq5k$. Let
$\{G_i\}_{i\in I}$ be a composition series of $G$, where $G_i$ is a proper subgroups of $G_j$ for
$i<j$. Let $S_1$, $S_2$, $\dots$, $S_{5k}$ be pairwise distinct nonabelian composition factors of
$G$. By Lemma~\ref{lfsg} every locally finite nonabelian simple group contains a finite insoluble
subgroup. Thus, we may choose finite subsets $X_1$, $X_2$, $\dots$, $X_{5k}$ of $G$ such that the
image of $X_i$ in $S_i$ generates an insoluble group. Suppose that $H$ is the finite subgroup of
$G$ generated by the union of the sets $X_1$, $X_2$, $\dots$, $X_{5k}$. Then $\{G_i\cap H\}_{i\in
I}$ is a subnormal series of $H$ having at least $5k$ insoluble factors. This contradicts
Proposition~\ref{FiniteCase}. The theorem is proved.

\noindent{\sl  Alexandr Buturlakin\\
Sobolev Institute of Mathematics,\\
Novosibirsk State University, \\
e-mail:buturlakin@math.nsc.ru}

~

\noindent{\sl Andrey V. Vasil$'$ev\\
Sobolev Institute of Mathematics,\\
Novosibirsk State University, \\
e-mail:vasand@math.nsc.ru}
\end{document}